\documentclass[a4paper,11pt,english]{article}

\usepackage[utf8]{inputenc}
\usepackage{amsthm}
\usepackage{amsmath,amssymb,amsfonts}
\usepackage{epsfig,color}
\usepackage{mathrsfs}
\usepackage{stmaryrd}
\usepackage{hyperref}
\usepackage{geometry}

\newcommand{\E}{\mathbb{E}}
\newcommand{\R}{\mathbb{R}}

\newtheorem{thm}{Theorem}

\title{A total variation version of Breuer--Major Central Limit Theorem under $\mathbb{D}^{1,2}$ assumption}

\author{J\"urgen Angst\footnote{Univ Rennes, CNRS, IRMAR - UMR 6625, F-35000 Rennes, France.\newline
\textcolor{white}{bla} \url{jurgen.angst@univ-rennes1.fr},\; \url{guillaume.poly@univ-rennes1.fr}}, \;
Federico Dalmao\footnote{DMEL, Cenur LN, Universidad de la República, Salto, Uruguay. \url{fdalmao@unorte.edu.uy} \newline 
\textcolor{white}{bla} This work was supported by the ANR grant UNIRANDOM, ANR-17-CE40-0008.} \; and Guillaume Poly$^*$}

\begin{document}

\maketitle
\begin{abstract}
In this note, we establish a qualitative total variation version of Breuer--Major Central Limit Theorem for a sequence of the type $\frac{1}{\sqrt{n}} \sum_{1\leq k \leq n} f(X_k)$, where $(X_k)_{k\ge 1}$ is a centered stationary Gaussian process, under the hypothesis that the function $f$ has Hermite rank $d \geq 1$ and belongs to the Malliavin space $\mathbb D^{1,2}$. This result in particular extends  the recent works of \cite{NNP21}, where a quantitative version of this result was obtained under the assumption that the function $f$ has Hermite rank $d= 2$ and belongs to the Malliavin space $\mathbb D^{1,4}$. We thus weaken the $\mathbb D^{1,4}$ integrability assumption to $\mathbb D^{1,2}$ and remove the restriction on the Hermite rank of the base function. While our method is still based on Malliavin calculus, we exploit a particular instance of Malliavin gradient called the sharp operator, which reduces the desired convergence in total variation to the convergence in distribution of a bidimensional  Breuer--Major type sequence.
\end{abstract}

\section{Framework and main result}\label{sec.intro}
Let us consider $X=(X_n)_{n \geq 1}$  a real-valued centered stationary Gaussian sequence with unit variance, defined on an abstract probability space $(\Omega, \mathscr{F}, \mathbb{P})$. Let $\rho : \mathbb N \to \mathbb R$ be the associated correlation function, in other words $\rho(|k-\ell|)=\mathbb E[ X_k X_{\ell}]$, for all $k,\ell \geq1$. We will also classically denote by $\mathcal N(0,\sigma^2)$ the law of a centered normal variable with variance $\sigma^2$.
Set $\gamma(d x):=(2 \pi)^{-1 / 2} e^{-x^2 / 2} d x$ the standard Gaussian measure on the real line and $\gamma_d=\otimes_{k=1}^d \gamma$ its analogue in $\mathbb R^d$. We then denote by $(H_m)_{m \geq 0}$ the family of Hermite polynomials which are orthogonal with respect to $\gamma$, namely $H_0 \equiv 1$ and
$$
H_m(x):=(-1)^m e^{\frac{x^2}{2}} \frac{d^m}{d x^m} e^{-\frac{x^2}{2}}, \quad m \geq 1.
$$
We denote by $L^2(\mathbb R, \gamma)$ the space of square integrable real functions with respect  to the Gaussian measure. Recall that a real function $f \in L^2(\mathbb R, \gamma)$ is said to have Hermite rank $d \geq 0$ if it can be decomposed as a sum of the form 
\[
f(x)= \sum_{m=d}^{+\infty} c_m H_m(x), \quad c_d \neq 0.
\]
For integers $k,p\geq 1$, we further denote by $\mathbb D^{k,p}(\mathbb R, \gamma)$ the Malliavin--Sobolev space consisting of the completion of the family of polynomial functions $q : \mathbb R \to \mathbb R$ with respect to the norm
\[
||q||_{k,p} :=\left| \int_{\mathbb R} \left( |q(x)|^p + \sum_{\ell=1}^k |q^{(\ell)}(x)|^p\right) \gamma(dx) \right|^{1/p},
\] 
where $q^{(\ell)}$ is the $\ell$-th derivative of $q$. 
Given a real function $f$, let us finally set 
\[
S_n(f):=\frac{1}{\sqrt{n}} \sum_{k=1}^n f(X_k).
\]
In this framework, the celebrated Central Limit Theorem (CLT) by Breuer and Major gives sufficient conditions on $\rho$ and $f$ so that the sequence 
$S_n(f)$ satisfies a CLT.

\begin{thm}[Theorem 1 in \cite{BM83}]\label{thm.bm}
If the function $f$ belongs to  $L^2(\mathbb R, \gamma)$ with Hermite rank $d\geq 1$ and if $\rho \in \ell^d(\mathbb N)$, i.e.
$
\sum_{\mathbb N} |\rho(k)|^d <+\infty,
$
then the sequence $(S_n(f))_{n \geq 1}$ converges in distribution as $n$ goes to infinity to a normal distribution $\mathcal N(0,\sigma^2)$, where the limit variance is given by
$$
\sigma^2:=\sum_{m=d}^{\infty} m ! c_m^2 \sum_{k \in \mathbb{Z}} \rho(k)^m,
$$
with $c_m$ being the coefficients appearing in the Hermite expansion of $f$.
\end{thm}

Recently, under mild additional assumptions, a series of articles has reinforced the above convergence in distribution into a convergence in total variation, with polynomial quantitative bounds, see e.g. \cite{KN19,NPY19,NZ21,NNP21}. Recall that the total variation distance between the distributions of two real random variables $X$ and $Y$ is given by 
\[
d_{\mathrm{TV}}(X,Y):=\sup_{A \in \mathcal B(\mathbb R)} | \mathbb P(X \in A)-\mathbb P(Y \in A)|,
\]
where the supremum runs over $\mathcal B(\mathbb R)$, the Borel sigma field on the real line.
To the best of our knowledge, the best statement so far in this direction is the  following

\begin{thm}[Theorem 1.2 in \cite{NNP21}]
Assume that $f \in L^2(\mathbb{R}, \gamma)$ has Hermite rank $d=2$ and that it belongs to $\mathbb{D}^{1,4}(\mathbb{R}, \gamma)$. Suppose that $\rho \in \ell^d(\mathbb N)$ and that the variance $\sigma^2$ of Theorem \ref{thm.bm} is positive.  Then, there exists a constant $C>0$ independent of $n$ such that
$$
d_{\mathrm{TV}}\left( \frac{S_n(f)}{\sqrt{\text{var}(S_n(f))}}, \mathcal N(0,1)\right) \leq \frac{C}{\sqrt{n}} \left[ \left(\sum_{|k| \leq n}|\rho(k)|\right)^{\frac{1}{2}}+\left(\sum_{|k| \leq n}|\rho(k)|^{\frac{4}{3}}\right)^{\frac{3}{2}}\right].
$$
\end{thm}
The goal of this note is to establish that the convergence in total variation in fact holds as soon as the function $f$ is in the Malliavin--Sobolev space $\mathbb D^{1,2}(\mathbb R, \gamma)$ and has Hermite rank $d \geq 1$.

\begin{thm}\label{thm.D12}
Suppose that $f \in \mathbb{D}^{1,2}(\mathbb{R}, \gamma)$ has Hermite rank $d \geq 1$. Suppose moreover that $\rho \in \ell^d(\mathbb N)$ and that the variance $\sigma^2$ of Theorem \ref{thm.bm} is positive. Then, as $n$ goes to infinity
$$
d_{\mathrm{TV}}\left( \frac{S_n(f)}{\sqrt{\text{var}(S_n(f))}}, \mathcal N(0,1)\right) \xrightarrow[n \to +\infty]{} 0.
$$
\end{thm}
\par
\bigskip
Note that, for the sake of simplicity, we only consider here a real Gaussian sequence $(X_n)_{n \geq 1}$ and a real function $f$ but our method is robust and would yield, under similar covariance and rank assumptions, a convergence in total variation for a properly renomalized sequence of the type $\sum_{k=1}^n f(X_k^1, \ldots, X_k^d)$ associated with a sequence of Gaussian vectors $(X_n)_{n \geq 1}$ with values in $\mathbb R^d$ and a function $f$ in the corresponding Malliavin--Sobolev space $\mathbb D^{1,2}(\mathbb R^d,\gamma_d)$. 
\par
\bigskip

The detailed proof of Theorem \ref{thm.D12} is the object of the next section and the rest of the paper. 
Unsurprisingly, we use the Malliavin--Stein approach to establish the CLT in total variation. However, our approach differs from the other works mentioned above in that we make use of the so called ``sharp gradient'', whose definition and main properties are recalled in Section \ref{sec.sharp}. With this tool at hand and in view of using Malliavin--Stein equation to characterize the proximity to the normal distribution, we shall see that the convergence in total variation in fact reduces to two rather simple steps\\
$i)$  a two-dimensional version of the classical Breuer--Major CLT  (i.e. in distribution not in total variation), see Section \ref{sec.2clt} ;\\
$ii)$ some elementary uniform integrability estimates, allowing to pass from a convergence in probability to a convergence in $L^1$, see Section \ref{sec.ui}.

\section{Proof of the main result}
As mentioned just above, the setting of the proof of Theorem \ref{thm.D12} is the one of Malliavin--Stein calculus. Note that for each fixed $n \geq 1$, the quantity of interest $S_n(f)$ involves only a finite number of Gaussian coefficients. So let us sketch the framework of Malliavin--Stein method in the finite dimensional setting, and we refer to \cite{Nualartbook} or \cite{NPbook} for a more general introduction. 
\subsection{A glimpse of Malliavin calculus}
Let us fix an integer $n\ge 1$ and let us place ourselves in the product probability space $(\mathbb{R}^n,\mathcal{B}(\mathbb{R}^n),\gamma_n)$ with $\gamma_n:=\otimes_{k=1}^n \gamma$, the $n$-dimensional standard Gaussian distribution on $\mathbb{R}^n$. Consider the classical \textit{Ornstein--Ulhenbeck} operator $\mathcal L_n:=\Delta-\vec{x}\cdot\nabla$ which is symmetric with respect to $\gamma_n$. We have then the standard decomposition of the $L^2-$space in Wiener chaoses, namely
\begin{eqnarray*}
L^2(\gamma_n)&=&\bigoplus_{k=0}^\infty \text{Ker}\left(\mathcal L_n+k\text{I}\right), \quad \text{with}\\
\text{Ker}\left(\mathcal L_n+k\text{I}\right)&=&\text{Vect}\left(\prod_{i=0}^n H_{k_i}(x_i)\Big{|}\sum_{i=0}^{n}k_i=n\right):=\hspace{-0.6 cm}\underbrace{\mathcal{W}_k.}_{k\text{-th Wiener chaos}}
\end{eqnarray*}
The square field or ``carr\'e du champ'' operator $\Gamma_n$ is then defined as the bilinear operator $\Gamma_n:=\left[\cdot,\cdot\right]=\nabla\cdot\nabla$. As a glimpse of the power of Malliavin--Stein approach in view of establishing total variation estimates, recall that if $F\in\text{Ker}\left(\mathcal L_n+k\text{I}\right)$ is such that $\mathbb{E}[F^2]=1$, then for some constant $C_k$ only depending on $k$, the total variation distance between the variable $F$ and a standard Gaussian can be upper bounded by
\[
d_{TV}\left(F,\mathcal{N}(0,1)\right)\le C_k \sqrt{\text{var}\left(\Gamma\left[F,F\right]\right)}.
\]
Via the notion of isonormal Gaussian process, the finite dimensional framework for Malliavin--Stein method sketched above can in fact be extended to the infinite dimensional setting giving rise to an Ornstein--Uhlenbeck operator $\mathcal L$ and an associated ``carr\'e du champ''  $\Gamma$, see e.g. Chapter 2 in \cite{NPbook}.

\subsection{The sharp gradient}\label{sec.sharp}

A detailed introduction to the sharp gradient can be found in Section 4.1 of the reference \cite{AP20}. We only recall here the basics which will be useful to our purpose.
Let us assume that $(N_k)_{k\ge 1}$ is an i.i.d. sequence of standard Gaussian variables on $(\Omega,\mathcal{F},\mathbb{P})$ which generate the first Wiener chaos. Without loss of generality, we shall assume that $\mathcal{F}=\sigma(N_k,\;k\ge 1)$. We will also need a copy $(\hat{\Omega},\hat{\mathcal{F}},\hat{\mathbb{P}})$ of this probability space as well as $(\hat{N}_i)_{i\ge 1}$ a corresponding i.i.d. sequence of standard Gaussian variables such that $\hat{\mathcal{F}}=\sigma(\hat{N}_k, \; k\ge 1)$. We will denote by $\hat{\mathbb E}$ the expectation with respect to the measure $\hat{\mathbb P}$.
For any integer $m\geq 1$ and any function $\Phi$ in the space $\mathcal{C}_b^1(\R^m,\R)$ of continuously differentiable functions with a bounded gradient, we then set
\begin{equation}\label{sharp operator}
{}^\sharp \Phi(N_1,\cdots,N_m):=\sum_{i=1}^m \partial_i \Phi(N_1,\cdots,N_m) \hat{N_i}.
\end{equation}
In Sections 4.1.1 and 4.1.2 of \cite{AP20}, it is shown that this \textit{gradient} is closable and extends to the Malliavin space $\mathbb{D}^{1,2}$, where 
\[
\mathbb D^{1,2} :=\left\lbrace F \in \mathbb L^{2}(\Omega, \mathcal F, \mathbb P), \, \mathbb E[ F^2] + \mathbb E\left[ ({}^\sharp F)^2  \right]<+\infty\right \rbrace.
\] 
The last space $\mathbb{D}^{1,2}$ is naturally the infinite dimensional version of the Malliavin--Sobolev space $\mathbb{D}^{1,2}(\mathbb R, \gamma)$ introduced in Section \ref{sec.intro} in the one-dimensional setting. In particular, Proposition 8 in the latter reference shows that
\[
\forall F\in\mathbb{D}^{1,2},\,\forall \phi \in \mathcal{C}^1_b(\R,\R)\,\,:\,\, {}^\sharp \phi(F)=\phi'(F) {}^\sharp F.
\]
Given $F\in\mathbb{D}^{1,2}$, taking first the expectation $\hat{\E}$ with respect $\hat{\mathbb P}$ and using Fubini inversion of sums yields the following key relation, for all $\xi \in \mathbb R$

\begin{equation}\label{Linksharpgamma}
\E\left(\exp\left(-\frac{\xi^2}{2}\Gamma[F,F]\right)\right)=\hat{\E}\E\left(\exp\left(i \xi {}^\sharp F \right)\right).
\end{equation}
By essence, via their Laplace/Fourier transforms, this key equation allows to relate the asymptotic behavior in distribution (or in probability if the limit is constant) of the carr\'e du champ $\Gamma[F,F]$ with the one of the sharp gradient ${}^\sharp F$. 
\par
\medskip
Finally, let us remark that by definition, the image $({}^\sharp X_k)_{k\ge 1}$ of our initial stationary sequence $(X_k)_{k\ge 1}$ by the sharp gradient is an independent copy of $(X_k)_{k\ge 1}$. We will write $({}^\sharp X_k)_{k\ge 1}=(\hat{X_k})_{k\ge 1}$ in the sequel.

\subsection{Convergence in probability via a two dimensional CLT}\label{sec.2clt}

Let us suppose that $f$ satisfies the assumptions of Theorem \ref{thm.D12}, namely $f \in \mathbb{D}^{1,2}(\mathbb{R}, \gamma)$ with Hermite rank $d \geq 1$, so that it can be decomposed as 
$f=\sum_{m=d}^\infty c_m H_m$ in $L^2(\mathbb R, \gamma)$. Let $\mathcal L^{-1}$ denote the pseudo-inverse of the Ornstein--Uhlenbeck operator and 
consider the pre-image
\[
g(x):=-\mathcal{L}^{-1}[f](x)=\sum_{m=d}^\infty \frac{c_m}{m} H_m(x).
\]
To simplify the expressions in the sequel, we set 
\[
F_n:=S_n(f)=\frac{1}{\sqrt{n}} \sum_{k=1}^n f(X_k), \;\; \text{and} \;\; G_n:=S_n(g)=-\mathcal{L}^{-1} F_n=\frac{1}{\sqrt{n}}\sum_{k=1}^n g(X_k).
\]
Now, take $(s,t,\xi)\in\R^3$ and let us apply the above key relation \eqref{Linksharpgamma} with the random variable $t F_n +s G_n$, we get
\begin{equation}\label{eq.key}
\mathbb{E}\left[\exp\left(-\frac{\xi^2}{2} \Gamma[t F_n+s G_n,t F_n+s G_n]\right)\right]=\hat{\E}\E\left[\exp\left(i\xi \left(t \, {}^\sharp F_n+ s \, {}^\sharp G_n \right)\right)\right].
\end{equation}
On the one hand, by bilinearity of the carr\'e du champ operator, we have
\begin{equation}\label{eq.bili}
\Gamma[t F_n+s G_n,t F_n+s G_n]=t^2\Gamma[F_n,F_n]+s^2 \Gamma[G_n,G_n]+2 t s \Gamma[F_n,-\mathcal{L}^{-1} F_n].
\end{equation}
On the other hand, the right hand side of Equation \eqref{eq.key} is simply the characteristic function under $\mathbb P \otimes \hat{\mathbb P}$ of  the couple $({}^\sharp F_n, {}^\sharp G_n)$ where
\[
\left( {}^\sharp F_n, \,{}^\sharp G_n\, \right) = \frac{1}{\sqrt{n}}\sum_{k=1}^n \left( f'( X_k) \hat{X_k}, g'( X_k) \hat{X_k} \right),
\]
is a ``Breuer--Major type" sequence with respect to the $\R^2-$valued centered stationary Gaussian process $(\hat{X_k},X_k)_{k\ge 1}$ and the $\R^2-$valued functional 
\[
(x,y)\mapsto \Psi(x,y):=(f'(x) y, g'(x) y).
\]
Since $f$ is in $\mathbb D^{1,2}(\mathbb R, \gamma)$, its derivative $f'$ is in $L^2(\mathbb R,\gamma)$ and $(\hat{X}_k)_{k\ge 1}$ and $(X_k)_{k\ge 1}$ are independent, therefore the  functional $\Psi$ is in $L^2(\mathbb R^2, \gamma_2)$ and the multivariate counterpart of the classical Breuer--Major Theorem applies, see Theorem 4 of \cite{arcones94}.

\bigskip

As a result, the bidimensional sequence $({}^\sharp F_n, {}^\sharp G_n)$ converges in distribution, under $\mathbb P \otimes \hat{\mathbb P}$, towards a bidimensional centered Gaussian vector with a symmetric semi-positive covariance matrix $\Sigma$. Therefore, from Equations \eqref{eq.key} and \eqref{eq.bili} and via the characterization of convergence in distribution in terms of Fourier transform, there exists real numbers $\lambda, \mu, \nu$ (depending on the limit covariance matrix $\Sigma$) such that for any $(s,t,\xi)\in\R^3$, as $n$ goes to infinity, we have
\[
\mathbb{E}\left[e^{-\frac{\xi^2t^2}{2}\Gamma[F_n,F_n]-\frac{\xi^2 s^2}{2} \Gamma[G_n,G_n]-\xi^2 t s \Gamma[F_n,-\mathcal{L}^{-1} F_n]}\right]
\xrightarrow[n\to\infty]~e^{-\frac{\xi^2}{2} \left(\lambda t^2+\mu s^2+2\nu t s\right)}.
\]
Since the above convergence is valid for any $\xi \in \mathbb R$, this shows in particular that for any fixed $(s,t) \in \mathbb R^2$, the sequence $\Gamma[t F_n+s G_n,t F_n+s G_n]$ converges in distribution (and thus in probability) towards the constant variable $\left(\lambda t^2+\mu s^2+2\nu t s\right)$. Choosing $s=t=1$, we thus get that $\Gamma[F_n+ G_n, F_n+ G_n]$ converges in probability towards $(\lambda+\mu+2\nu)$. Choosing $s=0$ and $t =1$, then $t=0$ and $s=1$, one deduce in the same manner that $\Gamma[F_n,F_n]$ and $\Gamma[G_n,G_n]$ both converge in probability towards $\lambda$ and $\mu$ respectively. Finally, by Equation \eqref{eq.bili}, one can conclude that the cross term 
\[
\Gamma[F_n, G_n]=\Gamma(F_n, -\mathcal{L}^{-1} F_n)=\hat{\mathbb E}\left[ {}^\sharp F_n {}^\sharp G_n \right]
\]
also converges in probability towards the constant limit variable $\nu$.

\subsection{Gaining some uniform integrability}\label{sec.ui}

Since our goal is to derive convergence in total variation of $F_n=S_n(f)$, the convergence in probability of the term $\Gamma[F_n,-\mathcal{L}^{-1} F_n]$ is not sufficient. Indeed, with Stein's Equation in mind, the lack of uniform integrability is a problem to deduce the following required asymptotic behavior for any $\phi\in\mathcal{C}^1_b(\mathbb{R})$, as $n$ goes to infinity
\[
\E\left[\phi'(F_n) \Gamma[F_n,-\mathcal{L}^{-1} F_n]\right]\approx \nu \, \E\left[\phi'(F_n)\right].
\]
In order to bypass this problem, let us go back to the two-dimensional classical Breuer--Major theorem associated with the functional $\Psi$ used in the last section. For any integer $p \geq 1$, let us denote by $\Psi_p$ the projection of $\Psi$ on the first $p-th$ chaoses. 
Applying Theorem 4 and Equation (2.43) of \cite{arcones94}, we get that there exists a constant $C>0$ (which depends only on the covariance structure of the underlying Gaussian process) such that 
\[
\sup_{n \ge 1} \E\hat{\mathbb E}\left[\left|\frac{1}{\sqrt{n}}\sum_{k=1}^n (\Psi-\Psi_p)(X_k,\hat{X}_k)\right]^2\right]\leq C \times \int_{\mathbb R^2} |(\Psi-\Psi_p)(x)|^2 \gamma_2(dx).
\]
Since $\Psi$ belongs to $L^2(\mathbb R^2, \gamma_2)$, the last term on the right hand side goes to zero as $p$ goes to infinity. As a result, uniformly in $n \geq 1$, the two-dimensional process
\[
\left( {}^\sharp F_n, \,{}^\sharp G_n\, \right) = \frac{1}{\sqrt{n}}\sum_{k=1}^n \Psi(X_k,\hat{X_k})
\]
can be approximated arbitrarily closely in $L^2(\mathbb P \otimes \hat{\mathbb P})$ by the following process which is finitely expanded on the Wiener chaoses
\[
Z_n^p:=(Z_n^{p,1},Z_n^{p,2}):=\frac{1}{\sqrt{n}}\sum_{k=1}^n \Psi_p(X_k,\hat{X_k}).
\]
Therefore, choosing $p\geq 1$ large enough, uniformly in $n\geq 1$, the product ${}^\sharp F_n \times {}^\sharp G_n$ can be approximated arbitrarily closely in $L^1(\mathbb P \otimes \hat{\mathbb P})$ by $\Delta_n^p := Z_n^{p,1} \times Z_n^{p,2}$. In other words, for any $\varepsilon>0$ and $p\geq 1$ large enough, we have 
\[
\sup_n \E\left[\left|\hat{\E}\left({}^\sharp F_n \times {}^\sharp G_n\right)-\hat{\E}\left(\Delta_n^p\right)\right|\right] \leq \sup_n\E\hat{\E}\left[\left|{}^\sharp F_n \times {}^\sharp G_n-\Delta_n^p\right|\right]<\varepsilon. 
\]
But mimicking the proof detailed in the previous Section \ref{sec.2clt} for the convergence in probability of $\Gamma[F_n, G_n]$ towards the constant variable $\nu$, one would then similarly get here that $\hat{\mathbb E}[\Delta_n^p]$ converges in probability under $\mathbb P$ towards a constant random variable $\nu_p \in \mathbb R$, and by construction $\lim_{p \to +\infty} \nu_p=\nu$. The crucial point here is that both random variables $\Delta_n^p$ and $\hat{\mathbb E}[\Delta_n^p]$ are now finitely expanded on the Wiener chaoses under $\mathbb{P}\otimes \hat{\mathbb{P}}$ and $\mathbb{P}$ respectively. Therefore, by hypercontractivity, the convergence in probability can be freely upgraded to the convergence in $L^q$ for every $q\ge 1$. In particular, as $n$ goes to infinity, the sequence $\hat{\mathbb E}[\Delta_n^p]$ converges in $L^1$ to the constant variable $\nu_p$.

\subsection{Conclusion}
We go back to Stein's Equation. Let $\phi\in\mathcal{C}^1_b(\mathbb{R})$ and $\varepsilon>0$.  Integrating by parts, for $p\geq 1$ large enough and by the results of the last section, we have 
\[
\begin{array}{ll}
\displaystyle{\left| \E\left[F_n\phi(F_n)\right]-\nu \, \E\left[\phi'(F_n)\right] \right|}  =\displaystyle{\left| \E\left[\phi'(F_n) \Gamma[F_n,-\mathcal{L}^{-1} F_n]\right] -\nu \, \E\left[\phi'(F_n)\right]\right|} \\
\\
 \displaystyle{=\left| \E\left[\phi'(F_n) \Gamma[F_n,G_n]\right] -\nu \, \E\left[\phi'(F_n)\right] \right|} \\
\\
 \displaystyle{= \left| \E\left[\phi'(F_n) \left(\Gamma[F_n,G_n]-\hat{\mathbb E}[\Delta_n^p]\right)\right]+\E\left[\phi'(F_n) \left(\hat{\mathbb E}[\Delta_n^p]-\nu_p\right)\right] + (\nu_p-\nu) \E\left[\phi'(F_n)\right]\right|}\\
 \\
  \displaystyle{\leq ||\phi'||_{\infty} \varepsilon +||\phi'||_{\infty}  \E\left[\left| \hat{\mathbb E}[\Delta_n^p]-\nu_p\right|\right] + ||\phi'||_{\infty}|\nu_p-\nu|.}
\end{array}
\]
As a result, letting first $n$ and then $p$ go to infinity, we get that uniformly in $\phi$ such that $||\phi'||_{\infty}\leq C$
\[
\limsup_{n \to +\infty} \left| \E\left[F_n\phi(F_n)\right]-\nu \, \E\left[\phi'(F_n)\right] \right| = 0.
\]
One can then classically conclude using Stein's approach for the convergence in total variation.


\end{document}